\documentclass[12pt,reqno]{amsart}
\usepackage{amsfonts,amssymb,amsmath}

\DeclareMathOperator{\bS}{\boldsymbol{\mathfrak{S}}}
\DeclareMathOperator{\bF}{\boldsymbol{\mathfrak{F}}}
\DeclareMathOperator{\bH}{\boldsymbol{\mathfrak{H}}}

\DeclareMathOperator{\size}{\mathit{d}}

\numberwithin{equation}{section} \numberwithin{theorem}{section}

\begin{document}

\title{Faces and Bases: Boolean Intervals}

\date{}
\author{Andrey O. Matveev}
\address{Data-Center Co., RU-620034, Ekaterinburg,
P.O.~Box~5, Russian~Federation} \email{aomatveev@\{dc.ru,
hotmail.com\}} 

\begin{abstract}
We consider redundant analogues of the $f$- and $h$-vectors of
simplicial complexes and present bases of $\mathbb{R}^{m+1}$
related to these ``long'' $f$- and $h$-vectors describing the face
systems $\Phi\subseteq\mathbf{2}^{\{1,\ldots,m\}}$; we list the
corresponding change of basis matrices. The representations of the
long $f$- and $h$-vectors of a face system with respect to various
bases are expressed based on partitions of the system into Boolean
intervals.
\end{abstract}

\maketitle

\section{Introduction and preliminaries}
\thispagestyle{empty}

Let $V$ be a finite set and let $\mathbf{2}^{V}$ denote the {\em
simplex\/} $\{F:\ F\subseteq V\}$. A family
$\Delta\subseteq\mathbf{2}^{V}$ is called an {\em abstract
simplicial complex\/} (or a {\em complex}) on the {\em vertex\/}
set $V$ if, given subsets $A$ and $B$ of $V$, the inclusions
$A\subseteq B\in\Delta$ imply $A\in\Delta$, and if
$\{v\}\in\Delta$, for any $v\in V$; see, e.g.,
\cite{BB,B,BH,BP,Hibi,MS,St1,Z}. If $\Gamma$ is a complex such
that $\Gamma\subset\Delta$ (that is, $\Gamma$ is a {\em
subcomplex\/} of $\Delta$) then the family $\Delta-\Gamma$ is
called a {\em relative simplicial complex},
see~\cite[\S{}III.7]{St1}.

If $\Psi$ is a relative complex then the sets $F\in\Psi$ are
called the {\em faces\/} of $\Psi$. The {\em dimension\/}
$\dim(F)$ of a face $F$ by definition equals $|F|-1$; the
cardinality $|F|$ is called the {\em size} of $F$. Let $\#$ denote
the number of sets in a family. If $\#\Psi>0$ then the {\em
size\/} $\size(\Psi)$ of $\Psi$ is defined by
$\size(\Psi):=\max_{F\in\Psi}|F|$, and the {\em dimension\/}
$\dim(\Psi)$ of $\Psi$ by definition is $\size(\Psi)-1$.

The row vector
$\pmb{f}(\Psi):=\bigl(f_0(\Psi),f_1(\Psi),\ldots,f_{\dim(\Psi)}\bigr)
\in\mathbb{N}^{\size(\Psi)}$, where $f_i(\Psi):=\#\{F\in\Psi:\
|F|=i+1\}$, is called the {\em $f$-vector\/} of $\Psi$. The row
{\em $h$-vector}
$\pmb{h}(\Psi):=\bigl(h_0(\Psi),h_1(\Psi),\ldots,h_{\size(\Psi)}\bigr)
\in\mathbb{Z}^{\size(\Psi)+1}$ of $\Psi$ is defined by
\begin{equation}
\sum_{i=0}^{\size(\Psi)}h_i(\Psi)\cdot\mathrm{y}^{\size(\Psi)-i}:=
\sum_{i=0}^{\size(\Psi)}f_{i-1}(\Psi)\cdot(\mathrm{y}-1)^{\size(\Psi)-i}\
.
\end{equation}

In this note we consider redundant analogues
$\pmb{f}(\Phi;|V|)\in\mathbb{N}^{|V|+1}$ and
$\pmb{h}(\Phi;|V|)\in\mathbb{Z}^{|V|+1}$ of the $f$- and
$h$-vectors that can be used in some situations for describing the
combinatorial properties of arbitrary {\em face systems}
$\Phi\subseteq\mathbf{2}^{V}$.

For a positive integer $m$, let $[m]$ denote the set
$\{1,2,\ldots,m\}$. We relate to a face system
$\Phi\subseteq\mathbf{2}^{[m]}$ the row vectors
\begin{align}
\label{eq:6}
\pmb{f}(\Phi;m):&=\bigl(f_0(\Phi;m),f_1(\Phi;m),\ldots,f_m(\Phi;m)\bigr)
\in\mathbb{N}^{m+1}\ ,\\ \label{eq:7}
\pmb{h}(\Phi;m):&=\bigl(h_0(\Phi;m),h_1(\Phi;m),\ldots,h_m(\Phi;m)\bigr)
\in\mathbb{Z}^{m+1}\ ,
\end{align}
where $f_i(\Phi;m):=\#\{F\in\Phi:\ |F|=i\}$, for $0\leq i\leq m$,
and the vector $\pmb{h}(\Phi;m)$ is defined by
\begin{equation}
\sum_{i=0}^m h_i(\Phi;m)\cdot\mathrm{y}^{m-i}:=\sum_{i=0}^m
f_i(\Phi;m)\cdot(\mathrm{y}-1)^{m-i}\ .
\end{equation}

Note that if $\Psi\subset\mathbf{2}^{[m]}$ is a relative complex
then we set $f_0(\Psi;m):=f_{-1}(\Psi):=\#\{F\in\Psi:\
|F|=0\}\in\{0,1\}$, $f_i(\Psi;m):=f_{i-1}(\Psi)$, for $1\leq
i\leq\size(\Psi)$ and, finally, $f_i(\Psi;m):=0$, for
$\size(\Psi)+1\leq i\leq m$.

Vectors~(\ref{eq:6}) and~(\ref{eq:7}) go back to analogous
constructions that appear, e.g., in~\cite{McMSh,McMWa}. In some
situations, these ``long'' $f$- and $h$-vectors either can be used
as an intermediate description of face systems or they can
independently be involved in combinatorial problems and
computations, see, e.g.,~\cite{M2}. Since the maps
$\Phi\mapsto\pmb{f}(\Phi;m)$ and $\Phi\mapsto\pmb{h}(\Phi;m)$ from
the Boolean lattice $\mathcal{D}(m)$ of all face systems (ordered
by inclusion) to $\mathbb{Z}^{m+1}$ are {\em valuations\/} on
$\mathcal{D}(m)$, the long $f$- and $h$-vectors can also be used
in the study of decomposition problems; here, a basic construction
is a {\em Boolean interval}, that is, the family
$[A,C]:=\{B\in\mathbf{2}^{[m]}:\ A\subseteq B\subseteq C\}$, for
some faces $A\subseteq C\subseteq[m]$.

We consider the vectors $\pmb{f}(\Phi;m)$ and $\pmb{h}(\Phi;m)$ as
elements from the real Euclidean space $\mathbb{R}^{m+1}$ of row
vectors. We present several bases of $\mathbb{R}^{m+1}$ related to
face systems and list the corresponding change of basis matrices.

See, e.g.,~\cite[\S{}IV.4]{Aigner} on valuations,
\cite[Chapter~5]{MS} on Alexander duality,
\cite[\S{}VI.6]{BarvinokConv}, \cite[Chapter~5]{BR},
\cite[\S{}II.5]{BH}, \cite[\S\S{}1.2, 3.6, 8.6]{BP},
\cite[\S{}III.11]{Hibi}, \cite[\S{}5.1]{McMSh}, \cite[\S\S{}II.3,
II.6, III.6]{St1}, \cite[\S{}3.14]{St2}, \cite[\S{}8.3]{Z} on the
Dehn-Sommerville relations, and~\cite{HJ} on matrix analysis.

\section{Notation}

Throughout this note, $m$ means a positive integer; all vectors
are of dimension $(m+1)$, and all matrices are $(m+1)\times(m+1)$
matrices. The components of vectors as well as the rows and
columns of matrices are indexed starting with zero. For a vector
$\pmb{w}$, $\pmb{w}^{\top}$ denotes its transpose.

If $\Phi$ is a face system, $\#\Phi>0$, then its {\em size\/}
$\size(\Phi)$ is defined by $\size(\Phi):=\max_{F\in\Phi}|F|$.

We denote the empty set by $\hat{0}$, and we use the notation
$\emptyset$ to denote the family containing no sets. We have
$\#\emptyset=0$, $\#\{\hat{0}\}=1$, and
\begin{align*}
\pmb{f}(\emptyset;m)&=\pmb{h}(\emptyset;m)=(0,0,\ldots,0)\ ,\\
\pmb{f}(\{\hat{0}\};m)&=\pmb{h}(\mathbf{2}^{[m]};m)=(1,0,\ldots,0)\
.
\end{align*}

$\boldsymbol{\iota}(m):=(1,1,\ldots,1)$;
$\boldsymbol{\tau}(m):=(2^m,2^{m-1},\ldots,1)$.

$\mathbf{I}(m)$ is the {\em identity matrix}.

$\mathbf{U}(m)$ is the {\em backward identity matrix} whose
$(i,j)$th entry is the Kronecker delta $\delta_{i+j,m}$.

$\mathbf{T}(m)$ is the {\em forward shift matrix} whose $(i,j)$th
entry is $\delta_{j-i,1}$.

If $\boldsymbol{\mathfrak{B}}:=(\pmb{b}_0,\ldots,\pmb{b}_m)$ is a
basis of $\mathbb{R}^{m+1}$ then, given a vector
$\pmb{w}\in\mathbb{R}^{m+1}$, we denote by
$[\pmb{w}]_{\boldsymbol{\mathfrak{B}}}:=\bigl(
\kappa_0(\pmb{w},\boldsymbol{\mathfrak{B}}),\ldots,
\kappa_m(\pmb{w},\boldsymbol{\mathfrak{B}})\bigr)\in\mathbb{R}^{m+1}$
the $(m+1)$-tuple satisfying the equality
$\sum_{i=0}^m\kappa_i(\pmb{w},\boldsymbol{\mathfrak{B}})\cdot\pmb{b}_i=\pmb{w}$.

\section{The long $f$- and $h$-vectors}

We recall the properties of vectors~(\ref{eq:6}) and~(\ref{eq:7})
described in~\cite{M1}.

\begin{itemize}
\item[\rm(i)]
The maps $\Phi\mapsto\pmb{f}(\Phi;m)$ and
$\Phi\mapsto\pmb{h}(\Phi;m)$ are valuations
$\mathcal{D}(m)\to\mathbb{Z}^{m+1}$ on the Boolean lattice
$\mathcal{D}(m)$ of all face systems (ordered by inclusion)
contained in $\mathbf{2}^{[m]}$.

\item[\rm(ii)]
Let $\Psi\subseteq\mathbf{2}^{[m]}$ be a relative complex.
\begin{align}
h_l(\Psi)&=\sum_{k=0}^l\binom{m-\size(\Psi)-1+l-k}{l-k}h_k(\Psi;m)\
,\ \ \ 0\leq l\leq\size(\Psi)\ ;\\ h_l(\Psi;m)&=(-1)^l\sum_{k=0}^l
(-1)^k\binom{m-\size(\Psi)}{l-k}h_k(\Psi)\ ,\ \ \ 0\leq l\leq m\ .
\end{align}

\item[\rm(iii)]
Let $\Phi\subseteq\mathbf{2}^{[m]}$.
\begin{itemize}
\item[\rm(a)]
\begin{align}
h_l(\Phi;m)&=(-1)^l\sum_{k=0}^l(-1)^k\binom{m-k}{l-k}f_k(\Phi;m)\
,\\ f_l(\Phi;m)&=\sum_{k=0}^l\binom{m-k}{l-k}h_k(\Phi;m)\ ,\ \ \
0\leq l\leq m\ .
\end{align}

\item[\rm(b)]
\begin{align}
h_0(\Phi;m)&=f_0(\Phi;m)\ ,\\
h_1(\Phi;m)&=f_1(\Phi;m)-mf_0(\Phi;m)\ ,\\
h_m(\Phi;m)&=(-1)^m\sum_{k=0}^m(-1)^k f_k(\Phi;m)\ ,\\
\pmb{h}(\Phi;m)\cdot\boldsymbol{\iota}(m)^{\top}&=f_m(\Phi;m)\ .
\end{align}

\item[\rm(c)]
\begin{equation}
\pmb{h}(\Phi;m)\cdot\boldsymbol{\tau}(m)^{\top}
=\pmb{f}(\Phi;m)\cdot\boldsymbol{\iota}(m)^{\top}=\#\Phi\ .
\end{equation}

\item[\rm(d)]
Consider the face system
\begin{equation}
\Phi^{\star}:=\{[m]-F:\ F\in\mathbf{2}^{[m]},\ F\not\in\Phi\}
\end{equation}
``dual'' to the system $\Phi$.

\begin{align}
h_l(\Phi;m)+(-1)^l\sum_{k=l}^m\binom{k}{l}h_k(\Phi^{\star};m)&=\delta_{l,0}\
,\ \ \ 0\leq l\leq m\ ;\\
h_m(\Phi;m)&=(-1)^{m+1}h_m(\Phi^{\star};m)\ .
\end{align}

If $\Delta$ is a complex on the vertex set $[m]$ then the complex
$\Delta^{\star}$ is called its {\em Alexander dual}. If
$\#\Delta>0$ and $\#\Delta^{\star}>0$ then
\begin{align}
h_l(\Delta;m)&=0\ ,\ \ \ 1\leq l\leq m-\size(\Delta^{\star})-1\
,\\
h_{m-\size(\Delta^{\star})}(\Delta;m)&=-f_{\size(\Delta^{\star})}
(\Delta^{\star};m)\ .
\end{align}
\end{itemize}
\end{itemize}

\section{Bases, and change of basis matrices}

We relate to the simplex $\mathbf{2}^{[m]}$ three pairs of bases
of the space $\mathbb{R}^{m+1}$. Let $\{F_0,\ldots,
F_m\}\subset\mathbf{2}^{[m]}$ be a face system such that
$|F_k|=k$, for $0\leq k\leq m$; here, $F_0:=\hat{0}$ and
$F_m:=[m]$.

The first pair consists of the bases $\bigl(\
\pmb{f}(\{F_0\};m),\pmb{f}(\{F_1\};m),\ldots,$
$\pmb{f}(\{F_m\};m)\ \bigr)$ and $\bigl(\
\pmb{h}(\{F_0\};m),\pmb{h}(\{F_1\};m),\ldots,\pmb{h}(\{F_m\};m)\
\bigr)$.

The bases $\bigl(\
\pmb{f}([F_0,F_0];m),\pmb{f}([F_0,F_1];m),\ldots,\pmb{f}([F_0,F_m];m)\
\bigr)$ and $\bigl(\
\pmb{h}([F_0,F_0];m),\pmb{h}([F_0,F_1];m),\ldots,
\pmb{h}([F_0,F_m];m)\ \bigr)$ compose the second pair.

The third pair consists of the bases $\bigl(\
\pmb{f}([F_m,F_m];m),\pmb{f}([F_{m-1},F_m];m),$
$\ldots,\pmb{f}([F_0,F_m];$ $m)\ \bigr)$ and $\bigl(\
\pmb{h}([F_m,F_m];m),\pmb{h}([F_{m-1},F_m];m),\ldots,\pmb{h}([F_0,$
$F_m];m)\ \bigr)$:

\begin{itemize}
\item[1)]
We use the notation $\bS_{m}$ to denote the {\em standard basis\/}
$\bigl(\boldsymbol{\sigma}(i;m):\ 0\leq i\leq m\bigr)$ of
$\mathbb{R}^{m+1}$, where
\begin{equation}
\boldsymbol{\sigma}(i;m):=(1,0,\ldots,0)\cdot\mathbf{T}(m)^i\ .
\end{equation}

We define a basis
$\bH^{\bullet}_m:=\bigl(\boldsymbol{\vartheta}^{\bullet}(i;m):\
0\leq i\leq m\bigr)$ of $\mathbb{R}^{m+1}$, where
\begin{equation}
\boldsymbol{\vartheta}^{\bullet}(i;m):=
\bigl(\vartheta^{\bullet}_0(i;m),\vartheta^{\bullet}_1(i;m),\ldots,
\vartheta^{\bullet}_m(i;m)\bigr)\in\mathbb{Z}^{m+1}\ ,
\end{equation}
by
\begin{equation}
\vartheta^{\bullet}_j(i;m):=(-1)^{j-i}\tbinom{m-i}{j-i}\ ,\ \ \
0\leq j\leq m\ .
\end{equation}

\item[2)]
Bases $\bF^{\blacktriangle}_m:=
\bigl(\boldsymbol{\varphi}^{\blacktriangle}(i;m):\ 0\leq i\leq
m\bigr)$ and $\bH^{\blacktriangle}_m:=
\bigl(\boldsymbol{\vartheta}^{\blacktriangle}(i;m):\ 0\leq i\leq
m\bigr)$ of $\mathbb{R}^{m+1}$ are defined in the following way:
\begin{equation}
\boldsymbol{\varphi}^{\blacktriangle}(i;m):=\bigl(\varphi^{\blacktriangle}_0(i;m),
\varphi^{\blacktriangle}_1(i;m),\ldots,\varphi^{\blacktriangle}_m(i;m)\bigr)
\in\mathbb{N}^{m+1}\ ,
\end{equation}
where
\begin{equation}
\varphi^{\blacktriangle}_j(i;m):=\tbinom{i}{j}\ ,\ \ \ 0\leq j\leq
m\ ,
\end{equation}
and
\begin{equation}
\boldsymbol{\vartheta}^{\blacktriangle}(i;m):=
\bigl(\vartheta^{\blacktriangle}_0(i;m),
\vartheta^{\blacktriangle}_1(i;m),\ldots,\vartheta^{\blacktriangle}_m(i;m)\bigr)
\in\mathbb{Z}^{m+1}\ ,
\end{equation}
where
\begin{equation}
\vartheta^{\blacktriangle}_j(i;m):=(-1)^{j}\tbinom{m-i}{j}\ ,\ \ \
0\leq j\leq m\ .
\end{equation}

The notations $\boldsymbol{\varphi}(i;m)$ and
$\boldsymbol{\vartheta}(i;m)$ were used in~\cite{M1} instead of
$\boldsymbol{\varphi}^{\blacktriangle}(i;m)$ and
$\boldsymbol{\vartheta}^{\blacktriangle}(i;m)$, respectively.

\item[3)]
The third pair consists of bases $\bF^{\blacktriangledown}_m:=
\bigl(\boldsymbol{\varphi}^{\blacktriangledown}(i;m):\ 0\leq i\leq
m\bigr)$ and $\bH^{\blacktriangledown}_m:=
\bigl(\boldsymbol{\vartheta}^{\blacktriangledown}(i;m):\ 0\leq
i\leq m\bigr)$ of $\mathbb{R}^{m+1}$ defined as follows:
\begin{equation}
\boldsymbol{\varphi}^{\blacktriangledown}(i;m):=
\bigl(\varphi^{\blacktriangledown}_0(i;m),
\varphi^{\blacktriangledown}_1(i;m),\ldots,
\varphi^{\blacktriangledown}_m(i;m)\bigr)\in\mathbb{N}^{m+1}\ ,
\end{equation}
where
\begin{equation}
\varphi^{\blacktriangledown}_j(i;m):=\tbinom{i}{m-j}\ ,\ \ \ 0\leq
j\leq m\ ,
\end{equation}
and
\begin{equation}
\boldsymbol{\vartheta}^{\blacktriangledown}(i;m):=
\bigl(\vartheta^{\blacktriangledown}_0(i;m),
\vartheta^{\blacktriangledown}_1(i;m),\ldots,
\vartheta^{\blacktriangledown}_m(i;m)\bigr)\in\mathbb{Z}^{m+1}\ ,
\end{equation}
where
\begin{equation}
\vartheta^{\blacktriangledown}_j(i;m):=\delta_{m-i,j}\ ,\ \ \
0\leq j\leq m\ .
\end{equation}
Note that $\bH^{\blacktriangledown}_m$ is up to rearrangement the
standard basis $\bS_{m}$.
\end{itemize}

Let $\mathbf{S}(m)$ be the change of basis matrix from $\bS_m$ to
$\bH^{\bullet}_m$: {\footnotesize
\begin{equation}
\mathbf{S}(m):=\begin{pmatrix}
\boldsymbol{\vartheta}^{\bullet}(0;m)\\ \vdots\\
\boldsymbol{\vartheta}^{\bullet}(m;m)
\end{pmatrix}\ ;
\end{equation}
} the $(i,j)$th entry of the inverse matrix $\mathbf{S}(m)^{-1}$
is $\tbinom{m-i}{j-i}$.

For any $i\in\mathbb{N}$, $i\leq m$, we have
\begin{align}
\boldsymbol{\vartheta}^{\bullet}(i;m)&=
\boldsymbol{\sigma}(i;m)\cdot\mathbf{S}(m)\ ,\\
\boldsymbol{\vartheta}^{\blacktriangle}(i;m)&=
\boldsymbol{\varphi}^{\blacktriangle}(i;m)\cdot\mathbf{S}(m)\ ,\\
\boldsymbol{\vartheta}^{\blacktriangledown}(i;m)&=
\boldsymbol{\varphi}^{\blacktriangledown}(i;m)\cdot\mathbf{S}(m)\
.
\end{align}

For any face system $\Phi\subseteq\mathbf{2}^{[m]}$, we have
\begin{align}
\label{eq:5} \pmb{h}(\Phi;m)&= \pmb{f}(\Phi;m)\cdot\mathbf{S}(m)
=\sum_{l=0}^m
f_l(\Phi;m)\cdot\boldsymbol{\vartheta}^{\bullet}(l;m) \ ,\\
\pmb{f}(\Phi;m)&= \pmb{h}(\Phi;m)\cdot\mathbf{S}(m)^{-1}\ .
\end{align}

The change of basis matrices corresponding to the bases defined
above are collected in Table~\ref{table:1}.

\section{Representations of the long $f$- and $h$-vectors
with respect to some bases}

If $\Phi\subseteq\mathbf{2}^{[m]}$ then we by~(\ref{eq:5}) have
\begin{equation}
\pmb{f}(\Phi;m)=[\pmb{h}(\Phi;m)]_{\bH^{\bullet}_m}\ ,
\end{equation}
and several observations follow:
\begin{align}
[\pmb{h}(\Phi;m)]_{\bH^{\blacktriangle}_m}&=
[\pmb{f}(\Phi;m)]_{\bF^{\blacktriangle}_m}\ ;\\
[\pmb{h}(\Phi;m)]_{\bH^{\blacktriangledown}_m}&=
[\pmb{f}(\Phi;m)]_{\bF^{\blacktriangledown}_m}\\ \nonumber
&=\pmb{h}(\Phi;m)\cdot \mathbf{U}(m)\ ;\\
[\pmb{f}(\Phi;m)]_{\bH^{\blacktriangledown}_m}&=\pmb{f}(\Phi;m)\cdot
\mathbf{U}(m)\\ \nonumber
&=[\pmb{h}(\Phi;m)]_{\bH^{\bullet}_m}\cdot \mathbf{U}(m)\ .
\end{align}

\section{Partitions of face systems into Boolean intervals,
and the long $f$- and $h$-vectors}

If
\begin{equation}
\label{eq:9}
\Phi=[A_1,B_1]\dot\cup\cdots\dot\cup[A_{\theta},B_{\theta}]
\end{equation}
is a partition of a face system $\Phi\subseteq\mathbf{2}^{[m]}$,
$\#\Phi>0$, into Boolean intervals $[A_k,B_k]$, $1\leq
k\leq\theta$, then we call the collection $\mathsf{P}$ of positive
integers $\mathsf{p}_{ij}$ defined by
\begin{equation}
\mathsf{p}_{ij}:=\#\{[A_k,B_k]:\ |B_k-A_k|=i,\ |A_k|=j\}>0
\end{equation}
the {\em profile} of partition~(\ref{eq:9}). If $\theta=\#\Phi$
then $\mathsf{p}_{0l}=f_l(\Phi;m)$ whenever $f_l(\Phi;m)>0$.
Table~\ref{table:2} collects the representations of the vectors
$\pmb{f}(\Phi;m)$ and $\pmb{h}(\Phi;m)$ with respect to various
bases.

\section{Appendix: Dehn-Sommerville type relations}

The $h$-vector of a complex $\Delta$ satisfies the {\em
Dehn-Sommerville relations\/} if it holds
\begin{equation}
h_l(\Delta)=h_{\size(\Delta)-l}(\Delta)\ ,\ \ \ 0\leq l\leq
\size(\Delta)
\end{equation}
or, equivalently (see, e.g.,~\cite[p.~171]{McMSh}),
\begin{equation}
h_l(\Delta;m)=(-1)^{m-\size(\Delta)}h_{m-l}(\Delta;m)\ ,\ \ \
0\leq l\leq m\ .
\end{equation}

We say, for brevity, that a face system
$\Phi\subset\mathbf{2}^{[m]}$ is a {\em DS-system} if the
Dehn-Sommerville type relations
\begin{equation}
\label{eq:2} h_l(\Phi;m)=(-1)^{m-\size(\Phi)}h_{m-l}(\Phi;m)\ ,\ \
\ 0\leq l\leq m
\end{equation}
hold. The systems $\emptyset$ and $\{\hat{0}\}$ are DS-systems.

If $\#\Phi>0$, then define the integer
\begin{equation}
\label{eq:3} \eta(\Phi):=\begin{cases}|\bigcup_{F\in\Phi}F|,
&\text{if $|\bigcup_{F\in\Phi}F|\equiv \size(\Phi)\pmod{2}$,}\\
|\bigcup_{F\in\Phi}F|+1, &\text{if
$|\bigcup_{F\in\Phi}F|\not\equiv \size(\Phi)\pmod{2}$.}
\end{cases}
\end{equation}
Note that, given a complex $\Delta$ with $v$ vertices, $v>0$, we
have
\begin{equation}
\eta(\Delta)=\begin{cases}v, &\text{if $v\equiv
\size(\Delta)\pmod{2}$,}\\ v+1, &\text{if $v\not\equiv
\size(\Delta)\pmod{2}$.}
\end{cases}
\end{equation}

Equality~(\ref{eq:2}) and definition~(\ref{eq:3}) lead to the
following observation: A face system $\Phi$ with $\#\Phi>0$ is a
DS-system if and only if for any $n\in\mathbb{P}$ such that
\begin{equation}
\label{eq:8}
\begin{split}
\eta(\Phi)&\leq n\ ,\\ n&\equiv\size(\Phi)\pmod{2}\ ,
\end{split}
\end{equation}
it holds
\begin{equation}
h_l(\Phi;n)=h_{n-l}(\Phi;n)\ ,\ \ \ 0\leq l\leq n\ ,
\end{equation}
or, equivalently,
\begin{equation}
\label{eq:11} \pmb{h}(\Phi;n)=\pmb{h}(\Phi;n)\cdot\mathbf{U}(n)\ ,
\end{equation}
that is, $\pmb{h}(\Phi;n)$ is a {\em left eigenvector\/} of the
$(n+1)\times(n+1)$ backward identity matrix corresponding to the
eigenvalue $1$.

We come to the following conclusion:

Let $\Phi$ be a DS-system with $\#\Phi>0$, and let $n$ be a
positive integer satisfying conditions~(\ref{eq:8}). Let
$l\in\mathbb{N}$, $l\leq n$.
\begin{itemize}
\item[\rm(i)]
\begin{align}
\nonumber \kappa_l\bigl(\pmb{h}(\Phi;n),
\bH^{\blacktriangle}_n\bigr)&=\kappa_l\bigl(\pmb{f}(\Phi;n),
\bF^{\blacktriangle}_n\bigr)\\ &=(-1)^{n-l}f_l(\Phi;n)\ ;\\
\nonumber \kappa_l\bigl(\pmb{h}(\Phi;n),
\bH^{\blacktriangledown}_n\bigr)&=\kappa_l\bigl(\pmb{f}(\Phi;n),
\bF^{\blacktriangledown}_n\bigr)\\ &=h_l(\Phi;n)=h_{n-l}(\Phi;n)\
;\\ \kappa_l\bigl(\pmb{h}(\Phi;n),
\bF^{\blacktriangle}_n\bigr)&=\kappa_l\bigl(\pmb{h}(\Phi;n),
\bF^{\blacktriangledown}_n\bigr)\ .
\end{align}
\item[\rm(ii)] If $\mathsf{P}$ is the profile of a partition of
$\Phi$ into Boolean intervals then the following equalities hold:
\begin{gather}
\sum_{i,j}{\mathsf p}_{ij}\cdot(-1)^{i+j}
\binom{j}{l-i}=(-1)^n\sum_{i,j}{\mathsf
p}_{ij}\cdot\binom{i}{l-j}\ ;\\ \sum_{i,j}{\mathsf
p}_{ij}\cdot(-1)^j\binom{n-i-j}{l-i}=(-1)^n\sum_{i,j}{\mathsf
p}_{ij}\cdot(-1)^j\binom{n-i-j}{l-j}\ ;\\ \nonumber
\sum_s\binom{s}{l}\sum_{i,j}{\mathsf p}_{ij}\cdot(-1)^j
\binom{n-i-j}{s-j}
\\ =(-1)^n\sum_s\binom{n-s}{l}\sum_{i,j}{\mathsf p}_{ij}\cdot(-1)^j
\binom{n-i-j}{s-j}\ .
\end{gather}
\end{itemize}

\newpage

{\footnotesize

\begin{table}[ht]
\caption{Change of basis matrices} \label{table:1}
\begin{center}
\begin{tabular}{|c|c|c|c|} \hline
\em Change of basis matrix & \em $(i,j)$th entry & \em Notation &
\em Case $m=3$\\ \hline\hline


from $\bS_m$ to $\bF^{\blacktriangle}_m$ & &
$\left(\begin{smallmatrix}\boldsymbol{\varphi}^{\blacktriangle}(0;m)\\
\vdots\\
\boldsymbol{\varphi}^{\blacktriangle}(m;m)\end{smallmatrix}\right)$
& \\

from $\bH^{\bullet}_m$ to $\bH^{\blacktriangle}_m$ &
$\binom{i}{j}$ & & $\left(\begin{smallmatrix}1&0&0&0\\ 1&1&0&0\\
1&2&1&0\\ 1&3&3&1\end{smallmatrix}\right)$ \\

from $\bH^{\blacktriangledown}_m$ to $\bF^{\blacktriangledown}_m$
& & &
\\ \hline


from $\bF^{\blacktriangle}_m$ to $\bS_m$ & &
$\left(\begin{smallmatrix}\boldsymbol{\varphi}^{\blacktriangle}(0;m)\\
\vdots\\
\boldsymbol{\varphi}^{\blacktriangle}(m;m)\end{smallmatrix}\right)^{-1}$
&
\\

from $\bH^{\blacktriangle}_m$ to $\bH^{\bullet}_m$ &
$(-1)^{i+j}\binom{i}{j}$ & & $\left(\begin{smallmatrix}1&0&0&0\\
-1&1&0&0\\ 1&-2&1&0\\ -1&3&-3&1\end{smallmatrix}\right)$
\\

from $\bF^{\blacktriangledown}_m$ to $\bH^{\blacktriangledown}_m$
& & &
\\ \hline\hline


from $\bS_m$ to $\bF^{\blacktriangledown}_m$ & &
$\left(\begin{smallmatrix}\boldsymbol{\varphi}^{\blacktriangledown}(0;m)\\
\vdots\\
\boldsymbol{\varphi}^{\blacktriangledown}(m;m)\end{smallmatrix}\right)$
& \\

from $\bH^{\bullet}_m$ to $\bH^{\blacktriangledown}_m$ &
$\binom{i}{m-j}$ & & $\left(\begin{smallmatrix}0&0&0&1\\ 0&0&1&1\\
0&1&2&1\\ 1&3&3&1\end{smallmatrix}\right)$  \\

from $\bH^{\blacktriangledown}_m$ to $\bF^{\blacktriangle}_m$ & &
&
\\ \hline


from $\bF^{\blacktriangledown}_m$ to $\bS_m$ & &
$\left(\begin{smallmatrix}\boldsymbol{\varphi}^{\blacktriangledown}(0;m)\\
\vdots\\
\boldsymbol{\varphi}^{\blacktriangledown}(m;m)\end{smallmatrix}\right)^{-1}$
& \\

from $\bH^{\blacktriangledown}_m$ to $\bH^{\bullet}_m$ &
$(-1)^{m-j-i}\binom{m-i}{j}$ & &
$\left(\begin{smallmatrix}-1&3&-3&1\\ 1&-2&1&0\\ -1&1&0&0\\
1&0&0&0\end{smallmatrix}\right)$  \\

from $\bF^{\blacktriangle}_m$ to $\bH^{\blacktriangledown}_m$ & &
& \\ \hline\hline


from $\bS_m$ to $\bH^{\blacktriangle}_m$ & $(-1)^j\binom{m-i}{j}$
&
$\left(\begin{smallmatrix}\boldsymbol{\vartheta}^{\blacktriangle}(0;m)\\
\vdots\\
\boldsymbol{\vartheta}^{\blacktriangle}(m;m)\end{smallmatrix}\right)$
& $\left(\begin{smallmatrix}1&-3&3&-1\\ 1&-2&1&0\\ 1&-1&0&0\\
1&0&0&0\end{smallmatrix}\right)$ \\ \hline


from $\bH^{\blacktriangle}_m$ to $\bS_m$ &
$(-1)^{m-j}\binom{i}{m-j}$ &
$\left(\begin{smallmatrix}\boldsymbol{\vartheta}^{\blacktriangle}(0;m)\\
\vdots\\
\boldsymbol{\vartheta}^{\blacktriangle}(m;m)\end{smallmatrix}\right)^{-1}$
& $\left(\begin{smallmatrix}0&0&0&1\\ 0&0&-1&1\\ 0&1&-2&1\\
-1&3&-3&1\end{smallmatrix}\right)$ \\ \hline\hline


from $\bS_m$ to $\bH^{\blacktriangledown}_m$ & &
$\left(\begin{smallmatrix}\boldsymbol{\vartheta}^{\blacktriangledown}(0;m)\\
\vdots\\
\boldsymbol{\vartheta}^{\blacktriangledown}(m;m)\end{smallmatrix}\right)$,
or &
\\

& $\delta_{m-i,j}$ & $\mathbf{U}(m)$, or &
$\left(\begin{smallmatrix}0&0&0&1\\ 0&0&1&0\\ 0&1&0&0\\
1&0&0&0\end{smallmatrix}\right)$
\\

from $\bH^{\blacktriangledown}_m$ to $\bS_m$ & &
$\left(\begin{smallmatrix}\boldsymbol{\vartheta}^{\blacktriangledown}(0;m)\\
\vdots\\
\boldsymbol{\vartheta}^{\blacktriangledown}(m;m)\end{smallmatrix}\right)^{-1}$
&
\\ \hline

\end{tabular}
\end{center}
\end{table}


$\quad$\newline $\quad$\newline $\quad$\newline $\quad$\newline
$\quad$\newline

\begin{center}
\begin{tabular}{|c|c|c|c|} \hline
\em Change of basis matrix & \em $(i,j)$th entry & \em Notation &
\em Case $m=3$\\ \hline\hline

from $\bS_m$ to $\bH^{\bullet}_m$ & $(-1)^{j-i}\binom{m-i}{j-i}$ &
$\left(\begin{smallmatrix}\boldsymbol{\vartheta}^{\bullet}(0;m)\\
\vdots\\
\boldsymbol{\vartheta}^{\bullet}(m;m)\end{smallmatrix}\right)$, or
& $\left(\begin{smallmatrix}1&-3&3&-1\\ 0&1&-2&1\\ 0&0&1&-1\\
0&0&0&1\end{smallmatrix}\right)$ \\ & & $\mathbf{S}(m)$&\\ \hline

from $\bH^{\bullet}_m$ to $\bS_m$ & $\binom{m-i}{j-i}$ &
$\left(\begin{smallmatrix}\boldsymbol{\vartheta}^{\bullet}(0;m)\\
\vdots\\
\boldsymbol{\vartheta}^{\bullet}(m;m)\end{smallmatrix}\right)^{-1}$,
or & $\left(\begin{smallmatrix}1&3&3&1\\ 0&1&2&1\\ 0&0&1&1\\
0&0&0&1\end{smallmatrix}\right)$\\ & & $\mathbf{S}(m)^{-1}$&\\
\hline\hline


from $\bF^{\blacktriangle}_m$ to $\bH^{\blacktriangle}_m$ &
$(-1)^j2^{m-j-i}\binom{m-i}{j}$ & &
$\left(\begin{smallmatrix}8&-12&6&-1\\ 4&-4&1&0\\ 2&-1&0&0\\
1&0&0&0\end{smallmatrix}\right)$ \\ \hline

from $\bH^{\blacktriangle}_m$ to $\bF^{\blacktriangle}_m$ &
$(-1)^{m-j}2^{i+j-m}\binom{i}{m-j}$ & &
$\left(\begin{smallmatrix}0&0&0&1\\ 0&0&-1&2\\ 0&1&-4&4\\
-1&6&-12&8\end{smallmatrix}\right)$\\ \hline\hline


from $\bF^{\blacktriangle}_m$ to $\bF^{\blacktriangledown}_m$ & &
& \\

from $\bF^{\blacktriangledown}_m$ to $\bF^{\blacktriangle}_m$ & &
& \\

& $(-1)^{m-j}\binom{m-i}{m-j}$ & &
$\left(\begin{smallmatrix}-1&3&-3&1\\ 0&1&-2&1\\ 0&0&-1&1\\
0&0&0&1\end{smallmatrix}\right)$ \\

from $\bH^{\blacktriangle}_m$ to $\bH^{\blacktriangledown}_m$ & &
& \\

from $\bH^{\blacktriangledown}_m$ to $\bH^{\blacktriangle}_m$ & &
& \\ \hline\hline


from $\bH^{\blacktriangle}_m$ to $\bF^{\blacktriangledown}_m$ &
$(-1)^{m-j}\sum_{s=\max\{m-i,m-j\}}^m\binom{i}{m-s}\binom{s}{m-j}$
& & $\left(\begin{smallmatrix}-1&3&-3&1\\ -1&4&-5&2\\ -1&5&-8&4\\
-1&6&-12&8\end{smallmatrix}\right)$\\ \hline

from $\bF^{\blacktriangledown}_m$ to $\bH^{\blacktriangle}_m$ &
$(-1)^{m-j}\sum_{s=0}^{\min\{m-i,m-j\}}\binom{m-i}{s}\binom{m-s}{j}$
& & $\left(\begin{smallmatrix}-8&12&-6&1\\ -4&8&-5&1\\ -2&5&-4&1\\
-1&3&-3&1\end{smallmatrix}\right)$  \\ \hline\hline


from $\bH^{\bullet}_m$ to $\bF^{\blacktriangle}_m$ &
$\sum_{s=0}^{\min\{i,j\}}\binom{i}{s}\binom{m-s}{m-j}$ & &
$\left(\begin{smallmatrix}1&3&3&1\\ 1&4&5&2\\ 1&5&8&4\\
1&6&12&8\end{smallmatrix}\right)$\\ \hline

from $\bF^{\blacktriangle}_m$ to $\bH^{\bullet}_m$ &
$(-1)^{i+j}\sum_{s=\max\{i,j\}}^m\binom{m-i}{m-s}\binom{s}{j}$ & &
$\left(\begin{smallmatrix}8&-12&6&-1\\ -4&8&-5&1\\ 2&-5&4&-1\\
-1&3&-3&1\end{smallmatrix}\right)$\\ \hline\hline


from $\bH^{\bullet}_m$ to $\bF^{\blacktriangledown}_m$ &
$2^{i+j-m}\binom{i}{m-j}$ & & $\left(\begin{smallmatrix}0&0&0&1\\
0&0&1&2\\ 0&1&4&4\\ 1&6&12&8\end{smallmatrix}\right)$\\ \hline

from $\bF^{\blacktriangledown}_m$ to $\bH^{\bullet}_m$ &
$(-2)^{m-j-i}\binom{m-i}{j}$ & &
$\left(\begin{smallmatrix}-8&12&-6&1\\ 4&-4&1&0\\ -2&1&0&0\\
1&0&0&0\end{smallmatrix}\right)$\\ \hline
\end{tabular}
\end{center}

}


{\footnotesize
\begin{table}[ht]
\caption{Representations (based on the profile of a partition of
$\Phi\subseteq\mathbf{2}^{[m]}$, $\#\Phi>0$, into Boolean
intervals) of $\pmb{f}(\Phi;m)$ and $\pmb{h}(\Phi;m)$ with respect
to various bases} \label{table:2}
\begin{center}
\begin{tabular}{|c|c|} \hline

\em $l$th component & \em Expression \\ \hline\hline

$f_l(\Phi;m)$ & $\sum_{i,j}{\mathsf p}_{ij}\cdot\binom{i}{l-j}$\\
\hline

$\kappa_l\bigl(\boldsymbol{\pmb{f}}(\Phi;m),
\bH^{\bullet}_m\bigr)$ & $\sum_s\binom{m-s}{m-l}\sum_{i,j}{\mathsf
p }_{ij}\cdot \binom{i}{s-j} $\\ \hline

$\kappa_l\bigl(\boldsymbol{\pmb{f}}(\Phi;m),
\bF^{\blacktriangle}_m\bigr)$ & $(-1)^l\sum_{i,j}{\mathsf p
}_{ij}\cdot(-1)^{i+j}\binom{j}{l-i}$\\ \hline

$\kappa_l\bigl(\boldsymbol{\pmb{f}}(\Phi;m),
\bH^{\blacktriangle}_m\bigr)$ &
$(-1)^{m-l}\sum_s\binom{s}{m-l}\sum_{i,j}{\mathsf
p}_{ij}\cdot\binom{i}{s-j} $\\ \hline

$\kappa_l\bigl(\boldsymbol{\pmb{f}}(\Phi;m),
\bF^{\blacktriangledown}_m\bigr)$ & $(-1)^{m-l}\sum_{i,j}{\mathsf
p}_{ij}\cdot(-1)^j\binom{m-i-j}{l-i}$\\ \hline

$\kappa_l\bigl(\boldsymbol{\pmb{f}}(\Phi;m),
\bH^{\blacktriangledown}_m\bigr)$ & $\sum_{i,j}{\mathsf p
}_{ij}\cdot\binom{i}{m-l-j}$\\ \hline \hline


$h_l(\Phi;m)$ & $(-1)^l\sum_{i,j}{\mathsf
p}_{ij}\cdot(-1)^j\binom{m-i-j}{l-j}$\\ \hline

$\kappa_l\bigl(\boldsymbol{\pmb{h}}(\Phi;m),
\bH^{\bullet}_m\bigr)$ & $\sum_{i,j}{\mathsf
p}_{ij}\cdot\binom{i}{l-j}$\\ \hline

$\kappa_l\bigl(\boldsymbol{\pmb{h}}(\Phi;m),
\bF^{\blacktriangle}_m\bigr)$ &
$(-1)^l\sum_s\binom{s}{l}\sum_{i,j}{\mathsf p}_{ij}\cdot(-1)^j
\binom{m-i-j}{s-j}$\\ \hline

$\kappa_l\bigl(\boldsymbol{\pmb{h}}(\Phi;m),
\bH^{\blacktriangle}_m\bigr)$ & $(-1)^l \sum_{i,j}{\mathsf
p}_{ij}\cdot(-1)^{i+j}\binom{j}{l-i}$\\ \hline

$\kappa_l\bigl(\boldsymbol{\pmb{h}}(\Phi;m),
\bF^{\blacktriangledown}_m\bigr)$ &
$(-1)^{m-l}\sum_s\binom{m-s}{l}\sum_{i,j}{\mathsf
p}_{ij}\cdot(-1)^j \binom{m-i-j}{s-j}$\\ \hline

$\kappa_l\bigl(\boldsymbol{\pmb{h}}(\Phi;m),
\bH^{\blacktriangledown}_m\bigr)$ & $(-1)^{m-l} \sum_{i,j}{\mathsf
p}_{ij}\cdot(-1)^j\binom{m-i-j}{l-i}$\\ \hline

\end{tabular}
\end{center}
\end{table}
}


\newpage

\begin{thebibliography}{9}
\bibitem{Aigner}
M.~Aigner, {\em Combinatorial Theory}, Grundlehren der
Mathematischen Wissenschaften, vol.~234, Springer-Verlag, Berlin,
1979.

\bibitem{BarvinokConv}
A.I.~Barvinok, {\em A Course in Convexity}, Graduate Studies in
Mathematics, vol.~54, American Mathematical Society, Providence,
RI, 2002.

\bibitem{BR}
M.~Beck and S.~Robins, {\em Computing the Continuous Discretely.
Integer-point Enumeration in Polyhedra}, Undergraduate Texts in
Mathematics, Springer-Verlag, {\em to appear}.

\bibitem{BB}
L.J.~Billera and A.~Bj\"{o}rner, {\em Face numbers of polytopes
and complexes}, in: {\em Handbook of Discrete and Computational
Geometry}, J.E.~Goodman and J.~O'Rourke (eds.) CRC Press, Boca
Raton, New~York, 1997, 291--310.

\bibitem{B}
A.~Bj\"{o}rner, {\em Topological methods}, in: {\em Handbook of
Combinatorics}, R.L.~Graham, M.~Gr\"{o}tschel and L.~Lov\'{a}sz
(eds.) {\em Vol.~2},  Elsevier, Amsterdam, 1995, 1819--1872.

\bibitem{BH}
W.~Bruns and J.~Herzog, {\em Cohen-Macaulay Rings, Second
edition}, Cambridge Studies in Advanced Mathematics, vol.~39,
Cambridge University Press, Cambridge, 1998.

\bibitem{BP}
V.M.~Buchstaber and T.E.~Panov,{\em Toricheskie Deistviya v
Topologii i Kombinatorike}. (in Russian) [{\em Torus Actions in
Topology and Combinatorics}] Moskovskii Tsentr Nepreryvnogo
Matematicheskogo Obrazovaniya, Moscow, 2004.

\bibitem{Hibi}
T.~Hibi, {\em Algebraic Combinatorics on Convex Polytopes},
Carslaw Publications, Glebe, Australia, 1992.

\bibitem{HJ}
R.A.~Horn and C.R.~Johnson, {\em Matrix Analysis}, Cambridge
University Press, Cambridge, 1986.

\bibitem{M1}
A.O.~Matveev, {\em Enumerating faces of complexes and valuations
on distributive lattices}, Discrete~Math.~Appl. {\bf 10} (2000)
no.~4, 403--421 (translation from Diskret.~Mat. {\bf 12} (2000)
no.~3, 76--94.)

\bibitem{M2}
A.O.~Matveev, {\em Faces and bases: Dehn-Sommerville type
relations}, preprint (2004).

\bibitem{McMSh}
P.~McMullen and G.C.~Shephard, {\em Convex Polytopes and the Upper
Bound Conjecture.} Prepared in collaboration with J.E.~Reeve and
A.A.~Ball. London Mathematical Society Lecture Note Series,
vol.~3, Cambridge University Press, London New~York, 1971.

\bibitem{McMWa}
P.~McMullen and D.W.~Walkup, {\em A generalized lower-bound
conjecture for simplicial polytopes}, Mathematika {\bf 18} (1971)
264--273.

\bibitem{MS}
E.~Miller and B.~Sturmfels, {\em Combinatorial Commutative
Algebra}, Graduate Texts in Mathematics, Springer-Verlag, 2004,
{\em to appear}.

\bibitem{St1}
R.P.~Stanley, {\em Combinatorics and Commutative Algebra, Second
edition}, Progress in Mathematics, Vol.~41, Birkhauser Boston,
Inc., Boston, MA, 1996.

\bibitem{St2}
R.P.~Stanley, {\em Enumerative Combinatorics, Vol. 1, Second
edition}, Cambridge Studies in Advanced Mathematics, vol.~49,
Cambridge University Press, Cambridge, 1997.

\bibitem{Z}
G.M.~Ziegler, {\em Lectures on Polytopes, Second edition},
Graduate Texts in Mathematics, vol.~152, Springer-Verlag,
New~York, 1998.
\end{thebibliography}
\end{document}